\documentclass[11 pt]{amsart}
 \usepackage{amsmath}
\usepackage{amssymb}
\usepackage{amsfonts}
     \usepackage{xcolor}
    \makeatletter
     \def\section{\@startsection{section}{1}%
     \z@{.7\linespacing\@plus\linespacing}{.5\linespacing}%
     {\bfseries \normalfont\scshape
     \centering
     }}
     \def\@secnumfont{\bfseries}
     \makeatother

   \newtheorem{theorem}{Theorem}[section]
\newtheorem{lemma}[theorem]{Lemma}

\newtheorem{proposition}[theorem]{Proposition}

\theoremstyle{definition}

\theoremstyle{remark}

\numberwithin{equation}{section}

\def\beq{\begin{equation}}
\def\eeq{\end{equation}}
  
\def\ben{\begin{eqnarray}}
\def\een{\end{eqnarray}}

\def \a{{\alpha}}
\def \b{{\beta}}
\def \D{{\Delta}}
 \def \d{{\delta}}

 \def \g{{\gamma}}

\def \l{{\lambda}}

 \def \m{{\mu}}
\def \s{{\sigma}}

\def \qq{{\qquad}}

\def\E{{\mathbb E}}

\def\R{{\mathbb R}}
\def \dd{{\rm d}}

at  6 pt
\scrollmode

 at 9,5 pt  \scrollmode

\font\sevenrm= cmr10 at 7,3 pt
 at 9,3 pt
 
\scrollmode

\def\ddate {\sevenrm \ifcase\month\or January\or
February\or March\or April\or May\or June\or July\or
August\or September\or October\or November\or December\fi\! {\the\day}, \!{\sevenrm\the\year}}
  
\begin{document}
   \title[A general decoupling inequality for finite Gaussian vectors]{ \rm { \bf   A general decoupling inequality for finite Gaussian vectors}
   } 
  
 \author{ Michel J.\,G.  WEBER}
 \address{Michel Weber: IRMA, UMR 7501, 10  rue du G\'en\'eral Zimmer, 67084
 Strasbourg Cedex, France}
  \email{michel.weber@math.unistra.fr \!}

      \begin{abstract}  
    We use Matrix Analysis to prove a general  decoupling inequality for finite  Gaussian vectors, in identifying a new region of   the  inherent $p$  exponent,  for the validity   of  this one.
 \end{abstract}
\maketitle


  \section{\bf Result}\label{s1}
     
 
  

  We are concerned in this work with finding optimal conditions for a    centered Gaussian vector  $X=\{X_i, 1\le i\le n\}$, to satisfy a decoupling inequality of type:     
    \beq\label{dec.ineq.intro}
 \E\Big( \prod_{i=1}^n f_i(X_i)\Big)
    \le    \mathcal Q(X,p)   \prod_{i=1}^n  \big\| f_i(X_i) \big\|_p, 
\eeq 
  for any complex-valued measurable  functions $f_1, \ldots, f_n$ such that $f_i\in L^p(\R)$, for all $1\le i\le n$, and where $p  $ is some real greater than 1, and $\mathcal Q(X,p) $ is an explicit constant.   
   \vskip 3 pt
  Call a  set of  reals $S\subset ]1,\infty)$ a $p$-admissible region, if for each $p\in  S$, Eq. \ref{dec.ineq.intro} holds with a corresponding explicit constant $\mathcal Q(X,p) $. 
  In the recent paper \cite{W3}
   the following decoupling inequality   is proved.  
   \begin{theorem}[\cite{W3},\,Th.\,2.3]\label{dec.ineq} Let $ X=\{X_i, 1\le i\le n\}$ be a centered Gaussian vector  with invertible  covariance matrix $C$, and let  $\s_i^2=\E X_i^2>0$ for each  $1\le i\le n$, $\underline{\g}=(\s^2_1, \ldots ,\s^2_n)$.   
Let $\b\ge 1$ be chosen so that $ \bar{\b}:=
\frac{(\max \s_i^2)}{(\min \s_i^2)}\vee\b>1$, and let $p$ be such that 
\begin{eqnarray}\label{cond}
p\, \ge \,   \bar{\b}\, p(X),\end{eqnarray}
where 
 $p(X)$ is the decoupling coefficient   of $X$   defined by
 $$ p(X)= \max_{i=1}^n\sum_{1\le j\le n} \frac{|\E X_iX_j|}{\E X_i^2}
\,.$$
 Then for any complex-valued measurable  functions $f_1, \ldots, f_n$ such that $f_i\in L^p(\R)$, for all $1\le i\le n$, the following inequality holds true, 
\begin{equation}\label{ineq}
\bigg|\,\E\Big( \prod_{i=1}^n f_i(X_i)\Big)\,\bigg|
\,\le \,\frac{ \big(\prod_{i=1}^n\s_i\big)^{\frac{1}{p}}}{\big(  1-  {1}/{\bar{\b}} \big)^{\frac{n}{2}(1-\frac{1}{p})}\det(C)^{\frac{1}{2p}}}\
  \prod_{i=1}^n  \big\| f_i(X_i) \big\|_p.
\end{equation}
 \end{theorem}
The $p$-admissible region is $S=[\bar{\b}\, p(X) ,\infty)$    with
$$\mathcal Q(X,p) =\frac{ \big(\prod_{i=1}^n\s_i\big)^{\frac{1}{p}}}{\big(  1-  {1}/{\bar{\b}} \big)^{\frac{n}{2}(1-\frac{1}{p})}\det(C)^{\frac{1}{2p}}}.$$ 
 The sharpness of the convenient condition Eq. \ref{cond} is however difficult to evaluate.   A cornerstone inequality in the proof is
\begin{equation}\label{det.lb.int.enonce}
\det\big(pI(\underline{\g})-C\big)  \,\ge\,  p^n\big(  1-\frac{1}{\bar{\b}} \big)^n\prod_{i=1}^n \s_i^2 ,
\end{equation} 
the notation  $ I(\underline{\g})$ being defined   below. See Section \ref{s4}.
  \vskip 6 pt
 In this paper we identify a general and new $p$-region of validity of   Eq. \ref{dec.ineq.intro}. This  region  is a   disconnected set and does not   express   in terms of    $C$,  but in terms of the columns of  another matrix directly  derived from $C^{1/2}$ and $ I(\underline{\g})$.  
 \vskip 2 pt We use  a theorem on simultaneous diagonalization. 
 The new region is more complicated than expected, but in the same time,   more intrinsically linked to $X$;  
  and we believe it is nearly optimal.  \vskip 2 pt     In comparison with Eq. \ref{det.lb.int.enonce} we obtain the exact formula,  
 \beq \label{kc.intro} \det\big(pI(\underline{\g})-C\big) =p^n \prod_{i=1}^n\sigma^2_i \Big(1-\frac{1}{\l_i}\Big).
\eeq  
  The reals $\l_j$ are positive and defined in Theorem \ref{dec.ineq.improv} below. 
\vskip 3 pt
{\it Notation}.  Let    $I_n$  be the $n\times n$ identity matrix and let $\underline b=(b_1, \ldots ,b_n)\in \R^n$. Let $  I(\underline b)$   denote the $n\times n$ diagonal matrix whose values on the diagonal are the corresponding values of $\underline b$. When $b_i\neq 0$ for each $i=1, \ldots, n$, we   use the notation $\underline b^{\a}=(b_1^{\a}, \ldots ,b_n^{\a})$, $\a$ real. Given a $n\times n$ matrix $M$,  ${}^tM$ denotes  the transpose  of $M$. We prove a general decoupling inequality under the minimal assumption that the covariance matrix $C$ of $X$  be invertible.  
\begin{theorem}
\label{dec.ineq.improv} Let $ X=\{X_i, 1\le i\le n\}$ be a centered Gaussian vector  with invertible  covariance matrix $C$, and let  $\s_i^2=\E X_i^2>0$ for each  $1\le i\le n$, $\underline{\g}=(\s^2_1, \ldots ,\s^2_n)$.  
 \vskip 2 pt {\rm (1)}
  Let  $1< p<\infty$.  There exists a  matrix $R$      expressed in terms of $C$ and $ I(\underline{\g})$ only,  with columns  $r^1, \ldots, r^n$  such that  letting 
\beq \label{det.B.cond.}
  \xi_j=\frac{\langle     I(\underline{\g})\,   r^j, r^j\rangle}{\langle C\, r^j, r^j\rangle}>0, \qq j=1, \ldots, n, 
\eeq
 \beq \label{det.B.cond.pos.sign.enonce}       S= \Big\{\frac{1}{\xi_i}, \ i=1,\ldots, n\Big\}^c\bigcap\Big\{p: \#\big\{i:\, 1\le i\le n;\,  p\xi_i< 1\big\}\  \hbox{\it  is even} \Big\} , \eeq
  if $p\in S$,  then 
        \beq\label{estimate.prod.c.enonce}
 \E\Big( \prod_{i=1}^n f_i(X_i)\Big)
    \le      \Big(\prod_{i=1}^n\s_i\Big)^{\frac{1}{p}}\det(C)^{-\frac{1}{2p}}   \Big(  \prod_{i=1}^n \Big|1-\frac{1}{p\, \xi_j }\Big|\Big)^{-\frac12(1-\frac{1}{p})} \
 \prod_{i=1}^n  \big\| f_i(X_i) \big\|_p,
\eeq 
 for any complex-valued measurable  functions $f_1, \ldots, f_n$ such that $f_i\in L^p(\R)$, for all $1\le i\le n$. 
\vskip 2 pt {\rm (2)}
 Condition \eqref{det.B.cond.pos.sign.enonce} is satisfied if $p> \max_{ j=1 }^n\, \frac{1}{\xi_i}$. 
 \vskip 3  pt {\rm (3)}
 Let $ \underline{\m}=\{\m_1,\ldots,\m_n\}$ be  the eigenvalues   $C$, and  $U$   an orthogonal matrix such that ${}^tUCU=I(\underline{\m})$. Let $D=I(\underline{\m}^{-1/2})$,  $W= UD$  so that  ${}^tW  C W=I$, and let matrix  $H= {}^tWKW$. The matrix  $H$ is real symmetric, positive definite.  
Let $\D$ be a real diagonal matrix, and  ${}^tVV=I$ such that ${}^tVHV=\D=I(\underline{\l})$.
  \vskip 2 pt Then $R$ equals the product\beq\label{comp.R} R=UDV,
  \eeq  $R$  is invertible, and we have the relations $\l_j=p\, \xi_j$, $j=1, \ldots, n$. 
 \end{theorem}
    Eq. \ref{estimate.prod.c.enonce}
 is    practical in Statistics,  by using   computing methods for computing eigenvalues from Numerical Analysis. See Section \ref{s5}.  \section{\bf Preparation}  \label{s2}
 
 \subsection{Brascamp-Lieb's inequality}The proposition below  follows from Theorem 6 in Brascamp and Lieb \cite{BL}. We also refer to \cite[Eq.\,4]{KLS} and   \cite[Prop.\,3.1]{W3}.

 \begin{proposition} \label{p1}Let $1\le p<\infty$. Let $A$ be a positive definite $n\times n$ matrix. Then for any measurable  functions $g_1, \ldots, g_n$ such that $g_i\ge 0$ and $g_i\in L^p(\R)$, $1\le i\le n$, the following inequality holds true, 
\begin{eqnarray}\label{BL.ineq}
\int_{\R^n}\Big( \prod_{i=1}^n g_i(x_i)\Big) \exp\Big\{ -\frac12\langle \underline x, A\underline x\rangle\Big\}\dd \underline x&\le & E_A\ \prod_{i=1}^n\Big( \int_{\R}    g_i(x)^p\dd x\Big)^\frac{1}{p},
\end{eqnarray}
where
\begin{eqnarray}\label{EB}
E_A&=&\sup_{b_i>0\atop i=1,\ldots,n}\frac{\int_{\R^n}\big( \prod_{i=1}^n \exp\{-\frac12 b_ix_i^2\}\big) \exp\{ -\frac12\langle \underline x, A\underline x\rangle\}\dd \underline x}{\prod_{i=1}^n\big( \int_{\R}   \exp\{-\frac{p}{2} b_ix^2\}\dd x\big)^\frac{1}{p}}\,.
\end{eqnarray}
\end{proposition}
  \subsection{Estimates of  $\boldsymbol{E_A}$}
  We have (\cite[p.\,705]{KLS},   or \cite[Remark 3.2]{W3}    for details),  
\begin{eqnarray*}\frac{\int_{\R^n}\big( \prod_{i=1}^n \exp\{-\frac12 b_ix_i^2\}\big) \exp\{ -\frac12\langle \underline x, A\underline x\rangle\}\dd \underline x}{\prod_{i=1}^n\big( \int_{\R}   \exp\{-\frac{p}{2} b_ix^2\}\dd x\big)^\frac{1}{p}} &=& \frac{\int_{\R^n} \exp\{ -\frac12\langle \underline x, (A+I(\underline b))\underline x\rangle\}\dd \underline x}{\prod_{i=1}^n\big( \frac{2\pi}{pb_i}  \big)^\frac{1}{2p}}
\cr &=&(2\pi)^{\frac{n}{2}(1-\frac{1}{p})}p^{\frac{n}{2p}} \,\frac{\prod_{i=1}^nb_i^\frac{1}{2p}}{ \det (A+I(\underline b))^\frac12}.
\end{eqnarray*}
So that
\begin{eqnarray}\label{EB1a}
E_A\,=\,(2\pi)^{\frac{n}{2}(1-\frac{1}{p})}p^{\frac{n}{2p}} \, \sup_{b_i>0\atop i=1,\ldots,n}\frac{\prod_{i=1}^nb_i^\frac{1}{2p}}{ \det (A+I(\underline b))^\frac12}.
\end{eqnarray}\vskip 2 pt 
\begin{proposition}[\cite{W3},\,Prop.\,3.4]\label{l1}
 \begin{eqnarray*}
E_A&\le&
\frac{(2\pi)^{\frac{n}{2}(1-\frac{1}{p})}}{\det(A)^{\frac12(1-\frac{1}{p})}} .
\end{eqnarray*}
\end{proposition}

\section{\bf Proof}\label{s3} (1)  By assumption  $C$ is real symmetric, positive definite. Thus    all eigenvalues   $\m_i$ of $C$ are positive, and there is an orthogonal matrix $U$ such that ${}^tUCU=I(\underline{\m})$. Also $C^{-1}= {}^tU I(\m^{-1})U$  is real symmetric. See Bellman  \cite[Th.\,2 p.\,54]{B}.   
   Set   
   \beq \label{B}B=C^{-1}-\frac{1}{p}I(\underline{\g}^{-1}) . 
   \eeq
 At first,
 \begin{align}\label{bl.a}
\int_{\R^n}\Big( \prod_{i=1}^n  |f_i(x_i)|  \Big)   \exp\Big\{\!\!&-\frac12\langle \underline x,  C^{-1}\underline x\rangle\Big\} \dd \underline x
\cr &=\int_{\R^n}\Big( \prod_{i=1}^n  |f_i(x_i)| e^{-x_i^2/(2p\s_i^2) })\Big)   \exp\Big\{\!\!-\frac12\langle \underline x, B\underline x\rangle\Big\} \dd \underline x .
\end{align}

By applying Proposition \ref{p1} with $A=B$, $g_i(x)= |f_i(x)| e^{-x^2/(2p\s_i^2) }$,   $i=1,\ldots, n$, it follows from Eq. \eqref{bl.a}  that 
  \begin{align*}
\int_{\R^n}\Big( \prod_{i=1}^n  |f_i(x_i)|  \Big)  \exp\Big\{\!\!-\frac12\langle \underline x, C^{-1}\underline x\rangle\Big\} \dd \underline x
   \le  \ E_{B} \,  \prod_{i=1}^n\Big( \int_{\R}    |f_i(x)|^p e^{-x^2/(2\s_i^2) } \dd x \Big)^\frac{1}{p},
\end{align*} 
if $B$ is positive definite. 
\vskip 3 pt 
Multiplying both sides by the factor $(2\pi)^{-\frac{n}{2} }$, next dividing each by $\det(C)^{\frac12}$ gives
\begin{align*}
\int_{\R^n}\Big( \prod_{i=1}^n  |f_i(x_i)| & \Big)   \exp\Big\{\!\!-\frac12\langle \underline x, C^{-1}\underline x\rangle\Big\}\frac{\dd \underline x}{(2\pi)^{\frac{n}{2}}\det(C)^{\frac12}} 
\cr \le &\ E_{B}\, \frac{ (2\pi)^{-\frac{n}{2}(1- \frac{1}{p} )}\big(\prod_{i=1}^n\s_i\big)^{\frac{1}{p}}}{ \det(C)^{\frac12}}  \prod_{i=1}^n\Big( \int_{\R}    |f_i(x)|^p e^{-x^2/(2\s_i^2) }\frac{\dd x}{\s_i\sqrt{2\pi}}\Big)^\frac{1}{p}.
\end{align*}
namely
\begin{eqnarray}\label{estimate.prod.a}
\E\Big( \prod_{i=1}^n f_i(X_i)\Big)\ \le \ E_{B}\, \frac{ (2\pi)^{-\frac{n}{2}(1- \frac{1}{p} )}\big(\prod_{i=1}^n\s_i\big)^{\frac{1}{p}}}{ \det(C)^{\frac12}}  \
 \prod_{i=1}^n  \big\| f_i(X_i) \big\|_p.
\end{eqnarray}  By  Proposition \ref{l1}    
 \begin{eqnarray}\label{EB.P34} 
E_{B}  \,\le\, 
\frac{(2\pi)^{\frac{n}{2}(1-\frac{1}{p})}}{\det(B)^{\frac12(1-\frac{1}{p})}}
, \end{eqnarray} 
if $B$ is positive definite, in which case
\begin{eqnarray}\label{estimate.prod.b}
 \E\Big( \prod_{i=1}^n f_i(X_i)\Big)
  &\le & \Big(\frac{(2\pi)^{\frac{n}{2}(1-\frac{1}{p})}}{\det(B)^{\frac12(1-\frac{1}{p})}}\Big)\, \frac{ (2\pi)^{-\frac{n}{2}(1- \frac{1}{p} )}\big(\prod_{i=1}^n\s_i\big)^{\frac{1}{p}}}{ \det(C)^{\frac12}}  \
 \prod_{i=1}^n  \big\| f_i(X_i) \big\|_p
 \cr &= &   \frac{  \big(\prod_{i=1}^n\s_i\big)^{\frac{1}{p}}}{  \det(B)^{\frac12(1-\frac{1}{p})}\det(C)^{\frac12}}  \
 \prod_{i=1}^n  \big\| f_i(X_i) \big\|_p.
\end{eqnarray}

Noticing that $ I(p \underline{\g})-C =C\big(C^{-1}I(p \underline{\g})-I\big)$ and $     C^{-1}I(p \underline{\g})-I=\big(C^{-1}- I(\frac{1}{p}\,\underline{\g}^{-1})\big)I(p \underline{\g})  $, we can write, 
\begin{eqnarray*}\det(I(p \underline{\g})-C)&=&\det\big( C (C^{-1}I(p \underline{\g})-I )\big)=\det(C)\det\big(   C^{-1}I(p \underline{\g})-I \big)
\cr &=& \det(C)\det\big(    C^{-1}- I(\frac{1}{p}\,\underline{\g}^{-1}))\det\big(I(p \underline{\g}) \big)
,
\end{eqnarray*} 
   and so we have
\begin{eqnarray}\label{bounddet}  \det(B)=  \det(   C^{-1}- I(\frac{1}{p}\,\underline{\g}^{-1})  ) &=&\frac{\det(I(p\,\underline{\g} )-C)}{\det(C)  \det(I(p\underline{\g} ))}\,=\,\frac{\det(I(p\,\underline{\g} )-C)}{p^n \det(C) \prod_{j=1}^n\s_j^2  }  .
\end{eqnarray} 
    \bigskip \par (2) By   Theorem 2 in   Franklin \cite[Sec.\,3.8]{Fr}, if 
$M$ and $K$ are two $n\times n$ real, symmetric, positive definite matrices, there exist  a   $n\times n$ real, invertible matrix $R$ and positive reals $\underline{\l}=(\l_1, \ldots ,\l_n)$    such that 
\beq \label{RMKR}{}^tRMR=I , \qq{}^tRKR=I(\underline{\l}) .  
\eeq  
Elementary properties of   determinants (products, transpose)  imply that 
$$\det( R )^2 =\frac{1}{\det(  M )} =\frac{\det(I(\underline{\l})  )}{ \det(  K )}  . $$
  Further
   $$I(\underline{\l})-I= {}^tRKR-{}^tRMR={}^tR(K-M)R.$$ Therefore 
   $$ \det(I(\underline{\l})-I )= \det( {}^tR(K-M)R ) =\det( R )^2\det( K-M ) .$$ It follows that
\beq \label{F.sec.38}\det( K-M )= \frac{ \det(I(\underline{\l})-I )}{\det( R )^2} =   \det(I(\underline{\l})-I )   \det(  M )=\frac{   \det(I(\underline{\l})-I ) \det(K)}{\det(I(\underline{\l})  )}   .
\eeq 


 
It is necessary in  what follows  to go to proof's details of Eq. \ref{RMKR}, namely the construction of $R$. More precisely, let $ \underline{\m}=\{\m_1,\ldots,\m_n\}$ be  the eigenvalues   $M$, and  $U$   an orthogonal matrix such that ${}^tUMU=I(\underline{\m})$. Let $D=I(\underline{\m}^{-1/2})$.   The effect of the transformation $L={}^tUKU$ is to transform $K$ into a matrix $L$ which is also real, symmetric, because ${}^tL=L$. Let $D=I(\underline{\m}^{-1/2})$. Then 
\beq\label{re1}{}^tD({}^tUCU)D=I, \qq {}^tD({}^tUKU))D= H.
\eeq 

\medskip \par 
The matrix $H$ being for the same reasons real, symmetric,  it  may be reduced  to a real diagonal matrix $\D$. If ${}^tVV=I$ and ${}^tVHV=\D=I(\underline{\l})$, Eq.\,\ref{re1} yields
\beq\label{re2}({}^tV{}^tD {}^tU)C(UDV)=I, \qq ({}^tV{}^tD {}^tU)K(UDV)= \D .
\eeq 
\vskip 3 pt Thus $R$ equals the product
$$ R=UDV,$$
  \ and $R$  is invertible. The $\l_j$ are computable: letting   $r^1, \ldots, r^n$ be  the columns of the matrix $R$, they equal the quotients
 \beq \label{comp.lj}\l_j=\frac{\langle K r^j, r^j\rangle}{\langle Mr^j, r^j\rangle}, \qq j=1, \ldots, n.   \eeq
See \cite[p.\,107]{Fr}. 
\bigskip \par (3) 
In the considered case  $K=pI(\underline{\g})$, $M =C$, we note that  $H$ takes the simpler form, computable as soon as $U$ is known, which belongs to the data, ($H={}^tD({}^tUKU))D$)
$$H=p\,I(\underline{\m}^{-1/2})({}^tUI(\underline{\g})U))I(\underline{\m}^{-1/2}) , $$  and $R=UI(\underline{\m}^{-1/2})V$,   $V$ being as above.
 We get 
\beq \label{comp.lj1} \l_j=\frac{\langle pI(\underline{\g}) r^j, r^j\rangle}{\langle Cr^j, r^j\rangle}=p\, \xi_j 
, \qq j=1, \ldots, n.   \eeq
By Eq. \ref{F.sec.38} and in view of Eq. \ref{det.B.cond.}, 
 we have the identity
\beq \label{kc} \det\big(pI(\underline{\g})-C\big) =p^n \prod_{i=1}^n\sigma^2_i \Big(1-\frac{1}{\l_i}\Big)=p^n \prod_{i=1}^n\sigma^2_i \Big(1-\frac{1}{p\, \xi_i}\Big)
.
\eeq 

Hence the lemma,
\begin{lemma} {\rm (1)} If \beq \label{det.B.cond}p \notin \Big\{\frac{1}{\xi_i}, \ i=1,\ldots, n\Big\},
 \eeq then
 $\det\big(pI(\underline{\g})-C\big) \neq 0$. 
\vskip 3 pt  {\rm (2)} Further if 
 \beq \label{det.B.cond.pos}p> \max_{ i=1,\ldots, n} \frac{1}{\xi_i}, 
 \eeq 
then  $\det\big(pI(\underline{\g})-C\big) > 0$.
\vskip 3 pt  {\rm (3)} Assume that   \eqref{det.B.cond} holds and that,  
 \beq \label{det.B.cond.pos.sign} \#\big\{i:\ 1\le i\le n;\  p\xi_i< 1\big\}\quad \hbox{\it  is even}, \eeq 
then  $\det\big(pI(\underline{\g})-C\big) > 0$. In particular (3) implies (2) since under assumption \eqref{det.B.cond.pos} the set in Eq. \ref{det.B.cond.pos.sign} has cardinality $0$. \end{lemma}
      In view of Eq. \ref{bounddet}, we deduce that under assumption \eqref{det.B.cond.pos}, we have
 $$\det(B)\,=\,\frac{\det(p\,I( \underline{\g} )-C)}{p^n(\prod_{j=1}^n\s_j^2)\det(C)   } \, =\, \frac{1 }{ \det(C) }  \prod_{i=1}^n \Big(1-\frac{1}{p\xi_i}\Big) > 0, $$
also
$$ \det(B)^{\frac12(1-\frac{1}{p})}\det(C)^{\frac12}=\Big(
\frac{1 }{ \det(C) }  \prod_{i=1}^n \Big(1-\frac{1}{p\xi_i}\Big)\Big)^{\frac12(1-\frac{1}{p})}\det(C)^{\frac12}.$$
$$=\Big(  \prod_{i=1}^n \Big(1-\frac{1}{p\xi_i}\Big)\Big)^{\frac12(1-\frac{1}{p})}\det(C)^{\frac{1}{2p}}. $$
Consequently by reporting in   Eq. \ref{estimate.prod.b},
 \begin{eqnarray}\label{estimate.prod.c}
 \E\Big( \prod_{i=1}^n f_i(X_i)\Big)
  &\le &  \frac{  \big(\prod_{i=1}^n\s_i\big)^{\frac{1}{p}}}{  \det(B)^{\frac12(1-\frac{1}{p})}\det(C)^{\frac12}}  \
 \prod_{i=1}^n  \big\| f_i(X_i) \big\|_p
 \cr &\le &    \Big(\prod_{i=1}^n\s_i\Big)^{\frac{1}{p}} \Big(  \prod_{i=1}^n \Big(1-\frac{1}{p\xi_i}\Big)\Big)^{-\frac12(1-\frac{1}{p})}\det(C)^{-\frac{1}{2p}}  \
 \prod_{i=1}^n  \big\| f_i(X_i) \big\|_p.
\end{eqnarray}

  \section{\bf Discussion}\label{s4}
  Estimate Eq. \ref{estimate.prod.b}   is a key step in the proof of Theorem \ref{dec.ineq}.   
    It remains in order to conclude to the result,   to bound from below $\det(B)$, which in view of Eq. \ref{bounddet}, amounts to bound from below $\det\big(p\,I(\underline{\g})-C\big)$. 
 This  is done in \cite{W3}  by   carefully applying Ostrowski's  lower bound (2) below:
 \begin{lemma} \label{l2}  Let $A=\{a_{i,j}, 1\le i,j\le n\}$.
{\rm (1)} Assume that
\begin{eqnarray}\label{hadamard.cond}
 |a_{i,i}|\,>\, \sum_{j=1\atop j\neq i}^n|a_{i,j}|, \qq \qq i=1,\ldots , n.
\end{eqnarray}
Then $\det(A)\neq 0$. 
\vskip 2 pt {\rm (2)} {\rm [\cite{Os},\,(2)]} Further under assumption \eqref{hadamard.cond},
\begin{eqnarray}\label{hadamard.lb}|\det(A)|\,\ge \,\prod_{i=1}^n \Big( |a_{i,i}|-\sum_{j=1\atop j\neq i}^n|a_{i,j}|\Big).
\end{eqnarray}\end{lemma} 
By assumption \eqref{cond},
 \begin{eqnarray*}
\sum_{1\le j\le n\atop j\neq i} |\E X_iX_j|\,\le \,\big(\frac{p}{ \bar{\b}}-1\big)\, \s_i^2, \qq \qq i=1,\ldots, n.
\end{eqnarray*}
 Letting $pI(\underline{\g})-C= \{d_{i,j},1\le i,j\le n\}$, we have $|d_{i,i}|=(p-1)\s_i^2$ and 
\begin{eqnarray*}
\sum_{j=1\atop j\neq i}^n|d_{i,j}|\,=\, \sum_{j=1\atop j\neq i}^n|\E X_iX_j|\,\le \, \big(\frac{p}{ \bar{\b}}-1\big)\, \s_i^2\,<\, |d_{i,i}|,\qq  i=1,\ldots , n,
\end{eqnarray*}
since $ \bar{\b}>1$. 
Thus
\begin{eqnarray*}|d_{i,i}|-\sum_{j=1\atop j\neq i}^n|d_{i,j}|&\ge&(p-1)\s_i^2-\Big(\frac{p}{ \bar{\b}} -1\Big)\, \s_i^2 
  \,=\,     p\s_i^2 \big(  1-\frac{1}{\bar{\b}} \big).
\end{eqnarray*}Applying inequality (2) yields in view of assumption \eqref{cond} on the choice of $p$ that   (\cite[Eq.\,(3.16)]{W3}) 
 \begin{equation}\label{det.lb.int.}
\det\big(pI(\underline{\g})-C\big)\,\ge\, \prod_{i=1}^n\big(p\s_i^2- \sum_{1\le j\le n\atop j\neq i} |\E X_iX_j|\big) \,\ge\,  p^n\big(  1-\frac{1}{\bar{\b}} \big)^n\prod_{i=1}^n \s_i^2 ,
\end{equation} 
 and finally by using \eqref{bounddet},   
 \begin{eqnarray*}\det(B)
&=&
\frac{\det\big(p\,I(\underline{\g})- C\big) }{p^n\,\det(C)\prod_{i=1}^n\s_i^2}
\cr &\ge &\frac{1 }{p^n\,\det(C)\prod_{i=1}^n\s_i^2}\,p^n\big(  1-\frac{1}{\bar{\b}} \big)^n\prod_{i=1}^n \s_i^2 
\,= \, \big(  1-\frac{1}{\bar{\b}} \big)^n\frac{1 }{\det(C)}.
\end{eqnarray*}

 \section{\bf Concluding Remarks} \label{s5} 
We begin with some remarks     concerning applications of Theorem \ref{dec.ineq.improv}.  Eq. \ref{kc} is a   simple case of simultaneous diagonalization.  In the cited theorem  \cite{Fr},  $M$ and $K$ are Hermitian, and only $M$ is supposed to be positive definite.  This is an important tool in Numerical Analysis. 
  The example of     mechanical systems of masses and springs is used to present the method, but also 
   serves 
as an application of it. See also \cite[Ch.\,4,\,Th.\,6]{B}. Applications of Theorem \ref{dec.ineq.improv} require to compute eigenvalues of $H$ (Eq. \ref{re2}).  There are many  computing methods for computing eigenvalues in Numerical Analysis, see \cite{B}. See \cite{Fr}, Sec. 7.15 for the   method of Francis and Kublanovskaya of computing eigenvalues, known as $QR$ method.    At this regard, an accurate way of reducing a $n\times n$ matrix $A$ to a right-triangular matrix  $R$ by 
performing a sequence of reflexions $U_{n-1}\ldots U_2 U_1 A=R$,  is described in the same remarkable guide-book \cite{Fr}, Sec. 7.14, and  is needed in the $QR$  method. The simple diagonalization is Gram-Schmidt process.   As demonstrated in Wilkinson's book \cite{Wi},  that process is 
 highly inaccurate in digital computation.      \medskip\par 
Our second remark is historical.  Claim (1) in  lemma \ref{l2} appeared in a book by Hadamard (1903), and  is in fact due  to L. L\'evy (1881) and in its general form to Desplanques (1887). We refer to Marcus and Minc \cite{MM},\, Ch.\,III, Sec.\,2 for more on its remarkable history. 
  \vskip 4 pt   For some classes of matrices,    a less known weaker condition  than \eqref{hadamard.cond} suffices to have $\det(A)\neq 0$.  Call an $n$-square matrix {\it irreducible} (or {\it indecomposable}) if it cannot be brought by a simultaneous row and column permutation to the form $\Big(\begin{matrix}X&Y\cr0&Z\end{matrix} \Big)$. If $A=\{a_{i,j}, 1\le i,j\le n\}$ is an irreducible matrix and 
\begin{eqnarray}\label{hadamard.cond.indec.}
 |a_{i,i}|\,\ge \, \sum_{j=1\atop j\neq i}^n|a_{i,j}|, \qq \qq i=1,\ldots , n,
\end{eqnarray}
with strict inequality for at least one $i$, then $\det(A)\neq 0$. This result is due to Taussky \cite{Ta}, see also  Marcus and Minc \cite[p.\,56]{MM} and Franklin \cite[pp.\,181-185]{Fr}. No   lower bound of $|\det(A)|$   seems known,   the way of proving   \eqref{hadamard.lb} failing in that case.  
\medskip\par


 \end{document}